\documentclass{amsart}

\usepackage{amsmath}
\usepackage{amssymb}
\usepackage{amsthm}
\usepackage{eufrak}
\usepackage{eucal}

\theoremstyle{plain}
\newtheorem{thm}{Theorem}[section]
\newtheorem{prp}[thm]{Proposition}

\newtheorem{fct}[thm]{Fact}
\newtheorem{cor}[thm]{Corollary}

\theoremstyle{definition}
\newtheorem{dfn}[thm]{Definition}

\newtheorem{ntn}[thm]{Notation}

\theoremstyle{remark}
\newtheorem{rmk}[thm]{Remark}

\newtheorem{qst}[thm]{Question}
\newtheorem{dsc}[thm]{Discussion}

\newtheorem{clm}{Claim}[thm]

\newtheorem{clm_asm}[clm]{Assumption}
\newtheorem{clm_rmk}[clm]{Remark}

\newtheorem{equ}[clm]{}
\newtheorem{tens_equ}[clm]{\m\bigotimes}

\theoremstyle{definition}

\newtheorem{imp_clm}[clm]{Claim}

\newenvironment{prf}{\begin{proof}}{\end{proof}}

\newcommand{\myqed}{\m{\textsc{q.e.d.}}}

\newcommand{\m}{\ensuremath}

\renewcommand{\l}{\m{\ell}}
\newcommand{\lam}{\m{\lambda}}
\newcommand{\om}{\m{\omega}}
\newcommand{\ka}{\m{\kappa}}
\newcommand{\al}{\m{\alpha}}
\newcommand{\be}{\m{\beta}}
\newcommand{\ga}{\m{\gamma}}
\newcommand{\rh}{\m{\rho}}
\newcommand{\ph}{\m{\varphi}}

\newcommand{\vrh}{\m{\varrho}}
\newcommand{\eps}{\m{\varepsilon}}

\newcommand{\then}{\m{\Rightarrow}}
\renewcommand{\iff}{\m{\Longleftrightarrow}}
\newcommand{\set}[1]{\ensuremath{\{#1\}}}
\newcommand{\seq}[1]{\ensuremath{\langle#1\rangle}}

\newcommand{\n}[1]{\ensuremath{\|#1\|}}

\newcommand{\ssup}[2]{\m{\sup_{#2}\set{#1}}}
\newcommand{\smax}[2]{\m{\max_{#2}\set{#1}}}

\newcommand{\func}[3]{\mbox{\ensuremath{#1:#2\to #3}}}
\newcommand{\sig}[2]{\mbox{\m{\Sigma_{#1}#2}}}

\newcommand{\supsub}[3]{#1^{#2}_{#3}}
\newcommand{\subsup}[3]{#1^{#3}_{#2}}

\newcommand{\dsubsup}[5]{#1^{#4^{#5}}_{#2_{#3}}}

\DeclareMathOperator{\tp}{tp}

\DeclareMathOperator{\cf}{cf}

\DeclareMathOperator{\len}{len}

\DeclareMathOperator{\tand}{\; and \;} \DeclareMathOperator{\tfor}{\; for \;}
 \DeclareMathOperator{\tst}{\; such \; that \;}
\DeclareMathOperator{\tif}{\; if \;} \DeclareMathOperator{\tin}{\; in \;}

\newcommand{\FD}{\mathfrak{D}}

\newcommand{\FK}{\m{\mathfrak{K}} \;}

\newcommand{\CL}{\m{\mathcal{L}} \;}

\newcommand{\CB}{\m{\mathcal{B} \;}}

\newcommand{\CG}{\m{\mathcal{G} \;}}

\newcommand{\BF}{\m{\mathbb{F} \;}}

\newcommand{\BR}{\m{\mathbb{R} \;}}
\newcommand{\BC}{\m{\mathbb{C} \;}}

\newcommand{\BZ}{\m{\mathbb{Z} \;}}

\newtheorem{mclm}[clm]{Main Claim}
\newenvironment{mclmprf}{\begin{proof}[Proof of the main claim]}{\end{proof}}

\renewcommand{\a}[1]{\m{\bar a_{#1}}}
\renewcommand{\b}[1]{\m {\bar b_{#1}}}
\newcommand{\dtpin}[3]{\mbox{\m{\tp(#1 #2, #3)}}}
\newcommand{\infseq}[1]{\mbox{\ensuremath{\langle #1_i : i < \omega \rangle}}}
\newcommand{\comut}[2]{\ensuremath{#1 #2 #1^{-1} = #2^{2}}}
\newcommand{\sop}[1]{\m{SOP_{#1}}}
\newcommand{\gen}[2]{\mbox{\m{\langle #1 \rangle_{#2}}}}

\newcommand{\bgen}[2]{\gen{\b{#1}}{#2}}

\newcommand{\dbgen}[3]{\gen{\b{#1}, \b{#2}}{#3}}
\newcommand{\Bm}[2]{\CB^{max}_{#1,#2}}
\newcommand{\hr}[1]{\m{\hat{r}_{#1}}}
\newcommand{\rt}[1]{\m{r'_{#1}}}

\newcommand{\famalg}[3]{\mbox{\m{#1 *_{#3} #2}}}
\newcommand{\dH}[2]{\m{H_{#1,#2}}}
\newcommand{\dbarphi}[2]{\mbox{\ensuremath{\varphi(\bar #1, \bar #2)}}}
\newcommand{\dbarphiun}[3]{\mbox{\ensuremath{\varphi_{#3}(\bar #1, \bar #2)}}}
\newcommand{\dphi}[2]{\mbox{\ensuremath{\varphi(#1, #2)}}}
\newcommand{\q}[2]{\mbox{\ensuremath{q(\bar #1, \bar #2)}}}
\newcommand{\p}[2]{\mbox{\ensuremath{p(\bar #1, \bar #2)}}}

\newcommand{\dalf}[2]{\m{f_{\al_#1, \al_#2}}}

\newcommand{\spc}[2]{\supsub{c}{#1}{#2}}
\newcommand{\norm}[2]{\mbox{\m{\parallel #1 \parallel_{#2}}}}
\newcommand{\sn}[1]{\mbox{\m{s(#1)}}}
\newcommand{\tterm}[4]{\mbox{\m{\tau_{#1,#2}(x_#3, y_{#4})}}}
\newcommand{\subtterm}[4]{\mbox{\m{\tau_{#1,#2}(#3, #4)}}}

\newcommand{\actson}%
           {\kern-1.36em{\rotatebox{-90}{\ensuremath{\circlearrowright}}}}
\title{Banach spaces and groups - order properties and universal models}
\author{Saharon Shelah \and
    Alex Usvyatsov}
\address{
Mathematics Department\\
Hebrew University of Jerusalem\\
91904 Givat Ram, Israel\\
}

\begin{document}
\maketitle

\begin{abstract}
    We deal with two natural examples of almost-elementary classes: the class of
    all Banach spaces (over \BR or \BC) and the class of all groups. We show both of
    these classes do not have the strict order property, and find the exact place of
    each one of them in Shelah's $SOP_n$ (strong order property of order $n$)
    hierarchy. Remembering the connection between this hierarchy and the existence of
    universal models, we conclude, for example, that there are ``few'' universal
    Banach spaces (under isometry) of regular cardinalities.

\end{abstract}


\section{Introduction and preliminaries}

    In this paper we deal with two very natural abstract elementary classes.
    An abstract elementary class (AEC) is a class of models for some dictionary
    (language) with a binary relation (order) defined on them, satisfying
    natural axioms saying that the order has a similar behavior to the
    "being an elementary submodel" relation in the first order case. One
    can see [Sh88] for definitions. AECs capture many examples of
    nonelementary classes, most of which being much more complicated
    to understand than the first order case, for example $\CL_{\ka,\om}$
    for $\ka > \om$. But our classes are not such good examples for seeing
    how general an AEC can be, as they are very similar to the first order
    case, and hardly even can be called ``abstract''. They have almost
    all the important properties that a regular elementary class has, and
    can be treated in a very similar way.

    The first class discussed in this paper is the class of all Banach
    spaces (real or complex), where ``being a submodel'' means just
    being a subspace - a linear subspace with the induced norm.
    As we have already mentioned, this class has very
    similar properties to an elementary class - it has amalgamation
    and the disjoint union property, locality of types (a type of an infinite
    sequence is determined by its finite subsequences) and it has
    compactness in a certain logic (positive strongly bounded formulae).
    For a special case of Banach spaces, positive strongly bounded formulae
    allow saying that the norm of a variable (or a constant) is in some \emph{compact}
    set of the reals, and they are closed under conjunction, disjunction
    and the existential quantifier. No negation or universal quantifiers
    are allowed. Itay Ben-Yaacov suggested in [BenYac] to call
    the classes having the above three properties CATs (compact
    abstract classes).

    Logic and model theory of Banach spaces was studied
    in detail by Henson and developed more by Henson and Iovino. The
    basic definitions can be found in [Iov], as well as
    a very good survey of the fundamentals of this
    theory, including Henson's compactness theorem for Banach spaces.
    Henson's logic allows more formulae (the ``admissible'' formulae
    are called ``positive bounded'') than the theory of CATs
    does, but there is also a price - the compactness theorem is ``local'',
    i.e. it is true inside every ball, but not in general (the norms
    of the variables have to be bounded by some uniform bound).

    The second class is even simpler to describe - it is the class of all
    groups, where being a submodel is just being a subgroup. Here we get
    all the properties of a CAT without even having to restrict our
    formulae. Compactness theorem holds trivially for groups in the
    regular first order language, as this is just a class of models of
    a universal first order theory. In fact, if we restrict ourselves
    to existentially closed groups, we will get a simple example of a
    Robinson theory (see [Hr]).

    One of main properties of a CAT, which makes the work really
    similar to the first order model theory, is the existence of
    a ``monster'' model, i.e. a model which is
    $\ka^*$-universal and $\ka^*$-homogeneous for $\ka^*$ much bigger
    than every cardinal mentioned in this paper
    (except $\ka^*$ itself, of course, which is also
    mentioned in the paper). Such models
    that are also very saturated and strongly homogeneous were called
    \emph{universal domains} by Hrushovski in [Hr], as they are the
    ``playground'' where all the work can be done. Every AEC with amalgamation
    and the disjoint union property has a ``monster''. For having
    a universal domain it is necessary to have locality of types and
    compactness, but this is also sufficient (see [Hr] or [BenYac]), so
    our classes have in fact very ``good'' monster models.

    So, what is the purpose of this paper? We've come to the questions
    asked and answered here from two directions that seem very different,
    though both deal with classification theory of models, and the
    connection between them was noticed by the first author
    in [Sh 500]. The original
    question arised from the project started by the first author of
    classifying elementary and non-elementary classes using the following
    ``test'' question:
    \begin{qst}
        In which regular $\lam \;$ can the theory/class in question have
        universal models, i.e. models that embed any other model of the
        same cardinality (and less)?
    \end{qst}
    This question was naturally asked about the class of Banach spaces.
    So the original question we were interested in is:
    \begin{qst}
    \label{qst_univ}
        In which regular $\lam \;$ can exist a universal Banach
        space?
    \end{qst}
    This question is certainly interesting for model theorists, as it
    has to do with classification theory of classes of models and asks
    how complicated a certain class is, but it
    is also of some interest to analysts researching Banach spaces themselves
    and their properties. Unlike saturated, big, compact models, etc,
    the concept and the importance of universal objects
    are well understood outside logic as well. Universal Banach spaces,
    for example, were studied by Banach himself in [Ban], and Szlenk
    in [Sz] showed there is no universal reflexive separable Banach
    space. Logicians and analysts do not always agree, though, on
    the question what a universal Banach space is. As the reader could
    have pointed out, an ``embedding'' of one model into another in our
    case is an isomorphism (in the usual sense in logic) of the first
    model onto a submodel (subspace) of the second one, i.e. a linear
    embedding which is also an \emph{isometry}. Analysts, on the other hand,
    usually allow more kinds of embeddings than
    just an isometry, and therefore right now the results presented here
    may not be of highest interest to them. Still, we will present a full
    answer to ~\ref{qst_univ}.

    The second question discussed in this paper arised from the joint
    interest of both authors in the $SOP_n$ (strong order property
    of order $n$) hierarchy for theories/classes defined by the first
    author in [Sh500]. We will recall the definitions:

    \begin{dfn}
    \label{dfn_sop}
    \begin{enumerate}
        \item
        We say that a formula \dbarphi{x}{y} exemplifies the
        \emph{strict order property (SOP)} in the model $M$ if it defines
        a partial order on $M$ with infinite indiscernible chains.

        \item
        We say that a formula \dbarphi{x}{y} exemplifies the
        \sop{n} (for $n \ge 3$) if it defines on $M$ a graph with
        infinite indiscernible chains and
        no cycles of size $n$.

        \item
        We say that a formula \dbarphi{x}{y} exemplifies the
        \sop{\le n} (for $n \ge 3$) if it defines on $M$ a graph with
        infinite indiscernible chains
        and no cycles of size smaller or equal to $n$.

        \item
        \label{dfn_AECop}
        We say an abstract elementary class \FK has
        \sop{}/\sop{n}/\sop{\le n} in a logic \CL
        if: there exists a formula \dbarphi{x}{y} in $\mathcal L$, such that
        for any infinite totally ordered set $I$ there is a model
        $M$ in \FK in which \sop{}/\sop{n}/\sop{\le n} is exemplified
        by \dbarphi{x}{y} with indiscernible chains of order type $I$.
    \end{enumerate}
    \end{dfn}

    \begin{rmk}
    \begin{enumerate}
        \item
            The idea of ~\ref{dfn_sop} (\ref{dfn_AECop}) is that
            without the compactness theorem, we have to demand
            indiscernible sequences of any length exemplifying an
            order property (which is not equivalent to just an
            infinite sequence). But we will see that in our cases
            this demand (and even the demand of indiscernibility)
            is unnecessary.

        \item
        Sometimes we will replace a formula in the above definitions
        by a type.

        \item
        One may view the \sop{n} hierarchy as finite approximations
        of the strict order property.

    \end{enumerate}
    \end{rmk}

    The \sop{n} hierarchy is connected in the following way to the
    well-known classes of stable and simple theories: any simple (and
    therefore any stable) theory/class does not have \sop{3}, which is
    the lowest property in the hierarchy.

    So a natural question can be:

    \begin{qst}
        Find natural examples of classes without the strict order
        property that are not trivial from the point of view of the
        \sop{n} hierarchy, i.e. are not simple.
    \end{qst}

    Here we show that both Banach spaces and groups provide good
    examples. Banach spaces have \sop{n} for every natural number $n$,
    and still do not have the $SOP$. Groups have the \sop{3}, but not the
    \sop{4}, which is even more surprising for such a complicated class.

    So what is the connection between the two questions? The answer
    was given by the first author in [Sh500]. Note that every AEC with
    amalgamation has a universal model in every regular $\lam \;$ satisfying
    $\lambda = \lambda^{<\lambda}$ or $\lambda = \mu^+$
    and $2^{< \mu} \le \lambda$. So it is definitely consistent that such
    a class has a universal model in every regular \lam - take $V=L$.
    Therefore ~\ref{qst_univ} can be asking only
    one thing - is it consistent that there is a universal Banach space
    of a regular cardinality $\lam \;$ that does \emph{not} satisfy any of the above
    equalities? In [Sh500] Shelah proves that the answer is negative for any
    class with \sop{4}. So finding the right location of the class of
    Banach spaces in the \sop{n} hierarchy gives automatically a full
    answer to the first question.

    An important question that has to be dealt with before we begin
    our discussion is - what should
    be the definition of the $SOP$ and \sop{n} for our abstract elementary
    classes? What do we mean by a ``formula'' that exemplifies an order
    property, what language do we allow? And is there a difference between
    demanding infinite chains, infinite indiscernible chains, or infinite
    indiscernible chains of any length?

    Considering the question of language - the case of groups
    is easy, here we deal with the regular first order formulae.
    What about Banach spaces? In fact, here we show
    the strongest possible results in each direction. We prove that
    \sop{n} is exemplified by using the most poor language - the language
    of CATs, i.e. positive \emph{strongly} bounded formulae, even
    quantifier free. In the other direction, we show that the strict order
    property can not be exemplified in \emph{any} ``locally'' compact
    language, i.e. can not be exemplified by any formula or type which
    is preserved under taking ultraproducts (where infinite elements
    are thrown away, and infinitesimal elements are divided out).
    In particular, Banach spaces do not have the strict order property
    in the rich Henson's language of positive bounded
    formulae, and even positive bounded types.

    What about the infinite chains that exemplify the order properties,
    do we have to demand explicitly existence of chains of any length, as
    this is usually done for an AEC,
    and is there a difference between regular and indiscernible chains?
    It turns out that in our cases, the compactness theorem solves all the
    problems. Suppose there exists in a ``monster'' model $M$
    an infinite sequence
    \infseq{\bar a} satisfying $i<j \implies M \models
    \dphi{\a{i}}{\a{j}}$. Then by compactness, there is such a sequence
    of any length (smaller than $\ka^*$) in $M$. Therefore, using the
    Erdos - Rado theorem, without loss of generality the original \om-sequence
    \infseq{\bar a} is also indiscernible (this trick is well-known
    for first order theories, and the proof works just the
    same for CATs - see [BenYac]). Now, again by
    compactness, there is an indiscernible chain as required of any
    length and order type (just a chain of the same type as
    \infseq{\bar a}). Note that one has to be very careful, as in the
    case of Banach spaces negation does not exist in the language, so
    not every technique can be used for finding indiscernible sequences
    (Ramsey's theorem won't help), but the process described above shows
    there is in fact no problem. So we can summarize:

    \begin{fct}
        Suppose there exists an infinite sequence \seq{\a{i}: i \in I}
        ($I$ - some infinite ordered set)
        in the universal domain $M$ of some compact abstract theory. Then
        for any infinite ordered set $J$,
        there exists an \emph{indiscernible} sequence \seq{\b{i}: i \in J}
        such that for all $n$ there exist $i_1 < \ldots < i_n$ in $I$
        satisfying $\tp(\b{0}, \ldots, \b{n-1}) =
        \tp(\dsubsup{\bar a}{i}{1}{}{}, \ldots, \dsubsup{\bar a}{i}{n}{}{})$.
        In particular, if \dphi{\a{i}}{\a{j}} for all $i<j \tin I$,
        the same thing holds for \seq{\b{i}: i \in J}.
    \end{fct}
    \begin{prf}
        Use Erdos-Rado and compactness,
        exactly like the case of $M$ a big model
        of a first order theory. For more details, see [BenYac].
    \end{prf}

    We will use the following immediate corollary:

    \begin{cor}
    \label{cor_indsop}
    \begin{enumerate}

        \item
        If $M$ is the universal domain of a compact abstract theory and we are
        interested in order properties exemplified in it,
        indiscernibility can be omitted from all the items of
        the definition ~\ref{dfn_sop}.

        \item
        If \FK is a compact abstract theory with the universal domain $M$
        (or just an abstract elementary class satisfying the compactness
        theorem for a logic \CL with the monster model $M$), \FK has
        \sop{}/\sop{n}/\sop{\le n} in \CL \iff it is exemplified in
        $M$ by some formula in \CL with indiscernible infinite chains
        of any order type (smaller than $\ka^*$) \iff it is exemplified
        in some $M_0 \in \FK$ (and therefore in $M$)
        by some formula in \CL with an infinite
        (not necessarily indiscernible) chain.

    \end{enumerate}
    \end{cor}


\section{Banach spaces}

Let \BF be either \BR or \BC.

\begin{ntn}
    We denote the ``monster'' Banach space (the universal domain
    of the compact abstract theory of Banach spaces) by \CB.
\end{ntn}

\begin{thm}
\label{thm_sopn}
    \CB has $SOP_n$ for all $n \ge 3$. Moreover, there is a positive
    strongly bounded quantifier free formula
    $\dbarphiun{x}{y}{n}$ exemplifying $SOP_{\le n}$
    in \CB with $\len{\bar x} = \len{\bar y} = 2$, such that
    $\dbarphiun{x}{y}{n+2} \vdash \dbarphiun{x}{y}{n}$.
\end{thm}
\begin{prf}
    Choose $n > 2$.

    First we define a seminormed space $B_0$. As a vector space over $\BF$, its
    basis is $\set{a_\al: \al < \om} \cup \set{b_\al : \al < \om}$.
    The seminorm is defined by $\sn{v} = \ssup{|f_\ga(v)|}{\ga < \om}$,
    where $f_\ga$ is a functional defined on the basis as follows:
    \[f_\ga(a_\al)=1 \tif \al < \ga, \; f_\ga(a_\al)=0 \tif \al \ge \ga \]
    \[f_\ga(b_\al)=0 \tif \al < \ga, \; f_\ga(b_\al)=1 \tif \al \ge \ga \]
    and extended to every $v \in B_0$ in the only possible way.
    Note that in fact $\sn{v} = \smax{|f_\ga(v)|}{\ga < \om}$ (i.e. the
    seminorm is finite). It is not a norm: it's easy to see that, for
    example, $\sn{a_0-a_1-b_1+b_0} = 0$. So we define $B_1$ as the normed
    space $B_1/\set{v: \sn{v}=0}$. Note that
    $\set{a_\al: \al < \om} \cup \set{b_\al : \al < \om}$ is no more
    a basis for $B_1$, though it certainly still is a set that generates
    the vector space. The following easy fact will be important for us:
    \begin{tens_equ}
    \label{equ_distseq}
        $\set{a_\al: \al < \om} \cup \set{b_\al : \al < \om}$ is
        a sequence of distinct non-zero elements in $B_1$.
    \end{tens_equ}
    Now denote
    the completion of $B_1$ by $B$ (of which we can think as of a
    subspace of the ``monster'' \CB).

    Now we define a \emph{term} in the language of Banach spaces
    (a positive bounded term) \tterm{n}{\l}{}{} by
    \[\tterm{n}{\l}{}{} = (n-2\l)x + (n-2\l+1)y \]
    Now define $c_{n,\l,\al} = \subtterm{n}{\l}{a_\al}{b_\al}$, i.e.
    \[c_{n,\l,\al} = (n-2\l)a_\al + (n-2\l+1)b_\al \]
    It is clear from the definitions that
    \[ |f_\ga(c_{n,\l,\al})| = n-2\l, \; \tif \al \le \ga \]
    and
    \[ |f_\ga(c_{n,\l,\al})| = n-2\l+1, \; \tif \al > \ga \].
    Therefore (check the calculation),
    \begin{tens_equ}
    \label{equ_1}
    \[
        \bigwedge_{\al < \be < \om}\bigwedge_{\l<n}
        \n{c_{n,\l+1,\be} - c_{n, \l, \al}} = 2
    \]
    \end{tens_equ}
    and
    \begin{tens_equ}
    \label{equ_2}
    \[
        \bigwedge_{\al < \be < \om}\bigwedge_{m \le n}
        \n{c_{n,m,\al} - c_{n, 0, \be}} = 2m+1
    \]
    \end{tens_equ}
    Now we define $\ph_n = \ph_n(x_1x_2, y_1y_2)$:
    \begin{eqnarray*}
        \ph_n =  \bigwedge_{\l<n}
        (\n{\tterm{n}{\l+1}{2}{2} - \tterm{n}{\l}{1}{1}} \le 2 )
        \; \& \\
          \bigwedge_{m \le n}
        (\n{\tterm{n}{m}{1}{1} - \tterm{n}{0}{2}{2}} \ge 2m+1 )
        \; \& \\
          \bigwedge_{m \le n}
        (\n{\tterm{n}{m}{1}{1} - \tterm{n}{0}{2}{2}} \le 2m+2 )
    \end{eqnarray*}
    \begin{clm_rmk}
        The last demand is not needed, as the reader will
        see in the proof, its only purpose is to make the
        formula strongly bounded. Readers who are interested
        only in Henson and Iovino's logic, can just omit it.
    \end{clm_rmk}

    Now we shall show that $\ph_n$ exemplifies $SOP_{\le n}$
    in $\CB$. First, by ~\ref{equ_1}, ~\ref{equ_2}, and of course
    ~\ref{equ_distseq},
    the sequence \seq{a_\al b_\al : \al<\om} verifies
    the first part of the definition (in $B$, of which we think
    as of a subspace  of $\CB$), i.e. it is an infinite chain
    of the graph defined by $\ph_n$ on \CB (by ~\ref{cor_indsop}, we
    don't need to prove indiscernibility).

    The only thing
    that is left to verify is that there are no cycles
    of length $\le n$ in this graph, and this is an immediate
    consequence of the triangle inequality (well hidden
    under the cover of long formulae):

    Suppose $2 < m \le n$, and suppose there are
    \seq{c_i, d_i : i \le m} in \CB such that
    $\CB \models \ph_n(c_i d_i, c_{i+1}d_{i+1}) \tfor i<m \tand
    \CB \models \ph_n(c_m d_m, c_0 d_0)$. Then in particular,
    from $\CB \models \ph_n(c_i d_i, c_{i+1}d_{i+1}) \tfor i<m$
    follows (taking only the ``first component'' of $\ph_n$):
    \begin{tens_equ}
    \label{equ_less}
    \[
        \bigwedge_{i<m}
        (\n{\subtterm{n}{i+1}{c_{i+1}}{d_{i+1}} -
        \subtterm{n}{i}{c_i}{d_i}} \le 2 )
    \]
    \end{tens_equ}
    On the other hand, $\CB \models \ph_n(c_m d_m, c_0 d_0)$ implies
    (taking only the ``second component'' of $\ph_n$):
    \begin{tens_equ}
    \label{equ_greater}
    \[
        \n{\subtterm{n}{m}{c_m}{d_m} -
        \subtterm{n}{0}{c_0}{d_0}} \ge 2m+1
    \]
    \end{tens_equ}
    But from \ref{equ_less} follows that
    $\n{\subtterm{n}{m}{c_m}{d_m} -
    \subtterm{n}{0}{c_0}{d_0}} \le
    \n{\subtterm{n}{m}{c_m}{d_m} -
    \subtterm{n}{m-1}{c_{m-1}}{d_{m-1}}} +
    \ldots +
    \n{\subtterm{n}{1}{c_1}{d_1} -
    \subtterm{n}{0}{c_0}{d_0}} \le 2m$,
    which contradicts \ref{equ_greater}.

    \begin{clm_rmk}
        Careful readers have probably pointed out that
        we actually showed that $\ph_n$
        exemplifies $SOP_{\le (n+1)}$ in \CB, but it
        doesn't matter for our discussion.
    \end{clm_rmk}

    We know now that $\ph_n$ exemplifies $SOP_{\le n}$.
    In order to complete the proof of the theorem, we need to
    show that $\ph_{n+2} \vdash \ph_{n}$ for all $n \ge 3$.
    For this, just note that \tterm{n}{\l}{}{} =
    \tterm{n+2}{\l+1}{}{}, and the rest follows immediately from
    the definition of $\ph_n$. \myqed
\end{prf}

The following corollary can be summarized as ``universal Banach spaces
in regular cardinals
exist only if they have to'', i.e. there are ``few'' universal Banach
spaces (under isometry).
\begin{cor}
    Suppose there exists a universal Banach space (under isometry)
    in $\lambda = \cf(\lambda)$.
    Then either $\lambda = \lambda^{<\lambda}$ or $\lambda = \mu^+$
    and $2^{< \mu} \le \lambda$.
\end{cor}
\begin{prf}
    $SOP_4$ is enough for this result - see [Sh500], Theorem 2.13.
\end{prf}

\begin{rmk}
    Note that the other direction of the last corollary is obvious -
    any abstract elementary class with amalgamation has a universal
    model in  every $\lam$ satisfying one of the above demands.
\end{rmk}

\begin{cor}
\label{cor:graph}
    There exists a positive strongly bounded quantifier free type
    type $p(\bar x, \bar y)$
    with $\len(\bar x) = \len(\bar y) = 2$, defining on \CB a graph
    with infinite (indiscernible) chains and no cycles at all.
\end{cor}
\begin{prf}
    Choose
    \[ \p{x}{y} = \bigwedge_{n \in \omega} \dbarphiun{x}{y}{2n+3} \]
    $p$ is consistent by compactness, as \dbarphiun{x}{y}{n+2}
    implies \dbarphiun{x}{y}{n}.
    Now, as  $\dbarphiun{x}{y}{n}$ exemplifies
    $SOP_{\le n}$ in \CB, and \seq{a_\al b_\al : \al < \om}
    from the proof of \ref{thm_sopn} is an
    infinite sequence ordered by $\ph_n$ for every $n$,
    the result is clear.
\end{prf}

\begin{dsc}
A natural question after we showed \ref{cor:graph} is: does \CB have the strict
order property? Or, a more general question: does having a (type-definable) graph as
in \ref{cor:graph} imply the strict order property (maybe also type-definable)?
Suppose we gave up compactness and allowed ourselves $L_{\omega_1, \omega}$
formulae, i.e. infinite disjunctions as well as infinite conjunctions. Then the
answer to the second question is certainly positive, as one can define the
transitive closure of a relation using an infinite disjunction, and the transitive
closure of \p{x}{y} is easily seen to be a partial order on \CB. But in our case the
implication is not clear, and in fact turns out to be false - we will give a
negative answer to the first question (and therefore to the second one). So the
compact abstract theory of Banach spaces turns out to be an interesting example of a
theory having a ``uniform'' definition of $SOP_n$, but yet without the $SOP$.
\end{dsc}

As nonstructure results for the class of Banach spaces are more
likely, the following one is rather surprising (and nice):

\begin{thm}
    \CB does not have the strict order property exemplified by a
    positive bounded type (in particular, \CB doesn't have the SOP
    exemplified by a p.b. formula).
\end{thm}
\begin{prf}
    Suppose towards a contradiction that \q{x}{y} is a ``compact''
    type which exemplifies \sop{} in \CB. So for every linear order $I$,
    there is an
    indiscernible sequence  $\seq{\a{i} : i \in I}$ which is linearly
    ordered by \q{x}{y}. We will choose $I=\BZ$.

    We denote $\len(\bar x) = \len(\bar y)$ in \q{x}{y} by $n$ and
    assume wlog that there exists $n^* < n$ such that
    $\bigwedge_{\ell<n^*} (a_{i,\ell} = a^*_\ell)$ for all $i \in I$
    and \seq{\bar a_{i,\ell} : n^* \le \ell < n, i \in I} is
    a linearly independent sequence. In other words, we assume

    \begin{clm_asm}\label{asm_ind}
        \seq{a^*_\ell: \ell<n^*} $\bigcup$
        \seq{\bar a_{i,\ell} : n^* \le \ell < n, i \in I}
        is a \emph{basis} for \gen{\a{i} : i \in I}{\CB}
    \end{clm_asm}

    Define for $k< \omega$, $B'_k = \gen{\a{k}, \a{k+1}}{\CB}$
    Denote for any $k_2 > k_1 + 1$, $B'_{k_1} \cap B'_{k_2}$ by $V^-$
    (generated by \seq{a^*_\ell : \l < n^*}).

    Pick $m< \omega$ and define $V_m$ as a vector subspace (over \BF)
    generated by \seq{\bar a_0, \bar  a_1, \ldots, \bar a_m} in \CB.
    Note that by
    ~\ref{asm_ind}, $V_m$ (as a vector
    space) is just
    a free amalgamation of $B'_0, \ldots, B'_{m-1}$ over
    $\gen{\a{1}}{}, \ldots, \gen{\a{k-1}}{}$ and $V^-$. We shall define
    three different norms on $V_m$. In order not to get confused between
    the original indiscernible sequence and the new normed space that
    we are going to define, we'll write \seq{\b{i} : i \le m} instead
    of \seq{\a{i} : i \le m}. Let $\func{h_1}{V_m}{\CB}$ and
    $\func{h_{-1}}{V_m}{\CB}$ be natural
    embeddings (isometries respecting the linear structure) such that
    for $0 \le k \le m$, $h_1(b_k) = a_k$ and $h_{-1}(b_k) = a_{-k}$.
    Let \func{g_{\l,k}}{\gen{\b{\l}}{}}{\gen{\b{k}}{}} be the natural
    isomorphism mapping \b{\l} onto \b{k} and $g_k = g_{0,k}$. Let
    \func{\subsup{h}{k}{\l}}{\gen{\b{k}}{}}{\gen{\a{\l}}{\CB}} be the
    natural isomorphism mapping \b{k} onto \a{\l}, i.e.
    $\subsup{h}{k}{\l} = h_1 \upharpoonright \gen{\b{\l}}{} \circ
    g_{k,\l}$.

    Now we define three different norms on $B'_k$ (for $k<m$).
    \norm{\cdot}{1} is a norm induced by $h_1$ (which is in fact
    the identity), \norm{\cdot}{-1} is induced by $h_{-1}$,
    \norm{\cdot}{0} is defined by $\max{\{ \norm{\cdot}{1},
    \norm{\cdot}{-1} \} }$. Now we expand these definitions to $V_m$:
    define for $t \in V_m$, $i \in \{ 0, 1, -1 \}$, $\norm{t}{i} =
    \inf{ \{ \sum_{k<m} \norm{t_k}{i} : t_k \in B'_k,
        \sum_{k<m} t_k = t \} }$.

    In fact, eventually we'll be interested only in \norm{\cdot}{1}.
    Our goal is to show that taking free amalgamations of
    $\gen{\a{0}, \a{1}}{\CB} \ldots \gen{\a{m-1}, \a{m}}{\CB}$ leads
    (in the limit - and here is where the compactness will be used)
    to a symmetric type. Two other norms are useful for
    showing the limit is symmetric, and their role will become
    clear in ~\ref{clm_conv}.

    Let $r', r'' \in \gen{\b{0}}{}$. Define for $0<k \le m$,
    $r_k = r' + g_k(r'')$. We will be interested in \norm{r_k}{i}
    for $i \in \set{0, 1, -1}$. Note that for $i \in \set{1, -1}$,
    by the definition of the norm \norm{\cdot}{i},
    for each $\epsilon > 0$,
    there are $t_p \in B'_p \tfor p<k$ such that
    $r_k = \sum_{p<k}t_p$ and
    $\norm{r_k}{i} + \epsilon
    \ge \sum_{p<k}\norm{t_p}{i} \ge \norm{r_k}{i}$.
    In the following claim we will assume that in fact one can
    find $t_p$ as above such that
    $\norm{r_k}{i} = \sum_{p<k}\norm{t_p}{i}$.

    \begin{clm}\label{clm_sum}
        Suppose $i \in \set{1,-1,0}$,
        $\norm{r_k}{i} = \sum_{p<k}\norm{t_p}{i}$
        where $t_p \in B'_p$ and $r_k = \sum_{p<k}t_p$.
        Then there exist $r'_p \in
        \gen{\b{p}}{}$ and $s_p \in V^- \tfor 0 \le p \le k \tst
         t_p = -r'_p + r'_{p+1} + s_p \tand  \tfor 0<p<k,
        r'_p \notin V^-$.
        Moreover, we may assume $r'_0 = -r'$ and
        $r'_k = g_k(r'')$, therefore $\sum_{p<k}{s_p} = 0$.
    \end{clm}

    \begin{prf}
        As $t_p \in B'_p$, we can write for every $p < k$,
        $t_p = \hat{r}_{p+1} - \check{r}_p + s_p$ for
        $\hat{r}_{p+1} \in \gen{\b{p+1}}{} \setminus V^-,
        \check{r}_p \in \gen{\b{p}}{} \setminus V^- \tand s_p \in V^-$.
        So we get $r_k = r' + g_k(r'') =
        - \check{r}_0 + \sum_{0<p<k}(\hat{r}_p - \check{r}_p) +
        \hat{r}_k + \sum_{p<k}s_p$. By ~\ref{asm_ind} and
        the definition of $V_m$,
        \seq{b^*_\ell: \ell<n^*} $\cup$
        \seq{\bar b_{i,\ell} : n^* \le \ell < n, i \le m} is
        a basis of $V_m$. As $\hat{r}_{p} \tand
        \check{r}_p$ are both elements of
        $\gen{\b{p}}{} \setminus V^-$, remembering the fact that
        $r_k = r' + g_k(r'')$, where $r' \in \gen{\b{0}}{} \tand
        r'' \in \gen{\b{k}}{}$, we get that necessarily
        $\hat{r}_{p} = \check{r}_p$. For  $0<p<k$, this
        is going to be $r'_p$. As the claim does not demand
        $r'_0, r'_k \notin V^-$, and we know that
        $r'+\check{r}_0 \in V^-$,
        as well as $g_k(r'')-\hat{r}_k$, by changing $s_0
        \tand s_{k-1}$, we may assume $\check{r}_0 = -r' \tand
        \hat{r}_k = g_k(r'')$. As
        \mbox{$r_k = -r'_0 + r'_k +\sum_{p<k}s_p$}, we get
        $\sum_{p<k}s_p = 0$, Q.E.D.
    \end{prf}

    Now we shall show

    \begin{clm}\label{clm_conv}
    \begin{enumerate}
    \item   For each $i \in \{0, 1, -1 \}$,
        $\norm{r_k}{i}$ is an ascending uniformly bounded sequence
        (the bound does not depend on $m$).

    \item   For each $j>1$, $m = j^2$,
        for each $i \in \{0, 1, -1 \}$,
        $\norm{r_m}{0} \ge \norm{r_m}{i} \ge
        (1+ \frac{2}{j})^{-1} \cdot \norm{r_j}{0}$

    \end{enumerate}
    \end{clm}

    \begin{prf}
    \begin{enumerate}
    \item   First we show the boundedness.
        $\norm{r' + g_k(r'')}{i} \le \norm{r'}{i} + \norm{g_k(r'')}{i}
        = \norm{h_1(r')}{\CB} + \norm{h_1(g_k(r''))}{\CB} =
        \norm{h_1(r')}{\CB} + \norm{h_1(r'')}{\CB}$. So as we see,
        the bound does not depend on $m$.

        Now suppose $k < \ell$. We aim to show that
        $\norm{r_k}{i} \le \norm{r_\ell}{i}$. First we'll prove this
        for $i=1$.

        As proving that
        for every  $\epsilon > 0$, $\norm{r_k}{1} \le
        \norm{r_\ell}{1} + \epsilon$ is enough,
        we may assume there exist $t_p \in B'_p$ such that
        $r_\ell = \sum_{p<\ell}t_p \tand
        \norm{r_\ell}{1} = \sum_{p<\ell}\norm{t_p}{1}$.
        Let $r'_p \in \gen{\b{p}}{}$ and
        $s_p \in V^- \tfor 0 \le p \le \ell$ be as in \ref{clm_sum}.

        Then
        \[
            r_k = r' + g_k(r'') = -r'_0 + g_k(r'')
        \]
        Therefore, by the definition of \norm{\cdot}{1} and $h_1$ being a linear function,
        \[
            \norm{r_k}{1} = \norm{-r'_0 + g_k(r'')}{1} =
            \norm{h_1(-r'_0 + g_k(r''))}{\CB} =
            \norm{h_1(-r'_0) + h_1(g_k(r''))}{\CB}
        \]
        But by indiscernibility of \a{i} in \CB and  the definitions of
        $g_k, g_\l, h_1$,
        \[
            \norm{h_1(-r'_0) + h_1(g_k(r''))}{\CB} =
            \norm{h_1(-r'_0) + h_1(g_\l(r''))}{\CB}
        \]
        So we get
        \[
        \begin{split}
            & \norm{r_k}{1} \le
            \norm{h_1(-r'_0) + h_1(g_\l(r''))}{\CB} =
            \norm{h_1(-r'_0 + r'_\l)}{\CB} =
            \norm{-r'_0 + r'_\l}{1} =
                \\
            & \norm{-r'_0 + r'_1 + s_0 - r'_1 + r'_2 + s_1 -
            \ldots - r'_{\l-1} + r'_\l + s_{\l-1} - \sum_{p<\l}s_p}{1} =
                \\
            & \norm{\sum_{p<\l}t_p -  \sum_{p<\l}s_p}{1}
        \end{split}
        \]
        Remembering that $\sum_{p<\l}{s_p} = 0$ (see ~\ref{clm_sum}),
        we conclude
        \[
            \norm{r_k}{1} \le
            \sum_{p<\l}\norm{t_p}{1} = \norm{r_\l}{1}
        \]
        finishing the proof for $i=1$.
        The same argument is used for $i=-1$, and the case $i=0$
        follows.

    \item
        Define (just for the proof) $\CB_{i,j} = \gen{\a{i},\a{j}}{\CB}$.
        Just as in case of \gen{\b{i},\b{j}}{}, we can define three
        norms on $\CB_{i,j} \;$: one is induced from the original norm
        on \CB (an analog of \norm{\cdot}{1}), the second one is
        induced from the norm on $\CB_{j,i}$, using the isomorphism
        from $\CB_{j,i}$ onto $\CB_{i,j}$ taking \a{i} onto \a{j} and
        vice versa (an analog of \norm{\cdot}{-1}). The third norm on
        $\CB_{i,j}$ (the one we will be actually interested in) will
        be denoted by \norm{\cdot}{\Bm{i}{j}}, and it is naturally
        an analog of \norm{\cdot}{0}, i.e. the maximum of the first two
        norms.

        So we start the proof with the following
        \begin{mclm}
        \label{main_clm}
            Suppose $m>k+1$ and $r = c_m - c_k \in
            \dbgen{k}{m}{}$, then $\norm{r}{1} \ge
            (1 + \frac{2}{m-k})^{-1}
            \cdot \norm{h_1(r)}{\Bm{k}{m}}$
        \end{mclm}

        \begin{mclmprf}

        First of all, wlog $k=0$.
        As in the previous proof, we assume
        the existence of $t_p \in B'_p \, \tfor p < m \tst
        r = \sig{}{t_p} \tand \norm{r}{1} = \sig{}{\norm{t_p}{1}}$.
        Therefore, by ~\ref{clm_sum}, there are $c_p \in
        \bgen{p}{} \setminus V^- \tfor p < m
        \tand s_p \in V^- \tfor p < m \tst$ for all $p$,
        $t_p = c_{p+1} - c_p + s_p$. Denote $\norm{t_p}{1} =
        \norm{h_1(t_p)}{\CB}$ by $\vrh_p$. So $\norm{r}{1} =
        \sig{p<m}{\vrh_p}$ and we aim to show
        \[
            \norm{h_1(r)}{\Bm{0}{m}} \le (1 + \frac{2}{m})\cdot
            \sig{}{\vrh_p}.
        \]
        Trivially (the triangle inequality)
        $\norm{h_1(r)}{\CB} \le \sig{}{\vrh_p} \le (1+\frac{2}{m})
        \cdot \sig{}{\vrh_p}$. Therefore it's left to show that
        \[
            \norm{h_{-1}(r)}{\CB} \le (1+\frac{2}{m})\cdot
            \sig{}{\vrh_p}
        \]

        Denote for $p<m \tand \al \in I$,
        $\subsup{c}{p}{\al} = \subsup{h}{p}{\al}(c_p)$.
        By the indiscernibility of \a{\al}, for all $\al < \be \in I$,
        \[
            \vrh_p = \norm{\subsup{c}{p+1}{\be} -
            \subsup{c}{p}{\al}}{\CB}
        \]
        Also, denote for some/all $\al < \be$
        \[
            \vrh^* = \norm{\subsup{c}{0}{\be} -
            \subsup{c}{m}{\al}}{\CB}
        \]

        For every $\al < \be \in I$ there is a functional
        $\func{f_{\al, \be}}{\CB}{\BF}$, such that
        \[
            \norm{f_{\al,\be}}{} = 1,
            f_{\al,\be}(\subsup{c}{0}{\be} -
            \subsup{c}{m}{\al}) = \vrh^*
        \]

        Choose \l \; such that $\vrh_\l$ is minimal. In particular,
        \begin{tens_equ}
        \label{tens_aver}
        \[
            \vrh_l \le \frac{1}{m}\sig{p<m}{\vrh_p}
        \]
        \end{tens_equ}

        Choose $\al_0<\al_1<\al_2<\al_3<\al_4$ in $I$.

        \[
        \begin{split}
            & \vrh^* = \n{\supsub{c}{\al_3}{0} -
            \supsub{c}{\al_1}{m}} =
            |{f_{\al_1, \al_3}(\supsub{c}{\al_3}{0} -
            \supsub{c}{\al_1}{m})}|
             = |\dalf{1}{3}(\spc{\al_3}{0}-\spc{\al_4}{1}) +
                \\
            & + \sum_{p=1}^{\l-1}
            \dalf{1}{3}(\spc{\al_4+p-1}{p} - \spc{\al_4+p}{p+1}) +
            \dalf{1}{3}(\spc{\al_4+\l-1}{\l} - \spc{\al_0-m+\l+1}{\l+1}) +
                \\
            & + \sum_{p=\l+1}^{m-2}
            \dalf{1}{3}(\spc{\al_0-m+p}{p} - \spc{\al_0-m+p+1}{p+1}) +
            \dalf{1}{3}(\spc{\al_0-m+(m-1)}{m-1} - \spc{\al_1}{m})|
            \le
                \\
            &\le \vrh_0 + \ldots + \vrh_{\l-1} +
            |\dalf{1}{3}(\spc{\al_4+\l-1}{\l} - \spc{\al_0-m+\l+1}{\l+1})| +
            \vrh_{\l+1} + \ldots + \vrh_{m-1}
        \end{split}
        \]
        The last inequality is true as
        $\norm{\dalf{1}{3}}{} = 1$. \\
        Denote $\be_1 = \al_4+\l-1$,
        $\be_2 = \al_0 - m + \l +1$. Find $\be_0 < \be_1 < \be_2 < \be_3$
        in $I$. \\
        Now note that
        \[
            \spc{\be_2}{\l} - \spc{\be_1}{\l+1} =
            (\spc{\be_2}{\l} - \spc{\be_3}{\l+1}) -
            (\spc{\be_0}{\l} - \spc{\be_3}{\l+1}) -
            (\spc{\be_0}{\l} - \spc{\be_1}{\l+1})
        \]
        Therefore,
        $\norm{\spc{\be_2}{\l} - \spc{\be_1}{\l+1}}{} \le 3\vrh_\l$.
        But (as $\norm{\dalf{1}{3}}{} = 1$),
        \[
            |\dalf{1}{3}(\spc{\al_4+\l-1}{\l} -
            \spc{\al_0-m+\l+1}{\l+1})| \le
            \norm{\spc{\be_2}{\l} - \spc{\be_1}{\l+1}}{}
            \le 3\vrh_\l
        \]
        Putting all the inequalities together (including ~\ref{tens_aver}), we
        conclude:
        \[
        \begin{split}
            & \vrh^* \le \vrh_0 + \ldots + \vrh_{\l-1} +
            |\dalf{1}{3}(\spc{\al_4+\l-1}{\l} -
            \spc{\al_0-m+\l+1}{\l+1})| +
            \vrh_{\l+1} + \ldots + \vrh_{m-1} \le
                \\
            & (\sum_{p<m}\vrh_p - \vrh_\l) + 3\vrh_\l
            = \sum_{p<m}\vrh_p + 2\vrh_\l \le
            \sum_{p<m}\vrh_p + 2 \cdot
            \frac{1}{m}\sum_{p<m}\vrh_p =
            (1+\frac{2}{m})\sum_{p<m}\vrh_p
        \end{split}
        \]
        which finishes the proof of the main claim.

        \end{mclmprf}

        Now assume $m=j^2 > 1$. We aim to show
        $\norm{r_j}{0} \le \norm{r_m}{1} \cdot (1+ \frac{2}{j})$.
        As usual, we assume $\norm{r_m}{1}
        = \sig{p<m}{\norm{t_p}{1}}$ for some $t_p \in B'_p$ satisfying
        $r = \sig{p<m}{t_p}$, where $t_p = -r'_p+r'_{p+1}+s_p$ as in
        ~\ref{clm_sum}.
        Denote for $\l \le j$, $\hr{\l}=
        g_{\l \cdot j \,, \l}(r'_{\l \cdot j})$, i.e. \hr{\l}
        is a copy of $r'_{\l \cdot j}$ in \gen{\b{\l}}{}.

        So by the definition of $r_j = r'+g_j(r'')$, we have
        \[
            \norm{r_j}{0} = \norm{r' + \hr{1} - \hr{1} +
            \hr{2} - \hr{2} + \ldots + \hr{j-1} - \hr{j-1}
            + g_j(r'')}{0}
        \]
        Now note that $g_j(r'') = \hr{j}$: $g_m(r'') = r'_m$
        (by ~\ref{clm_sum}),
        therefore $g_j(r'') = g_{m,j}(g_m(r'')) = g_{m,j}(r_m)$.
        Remembering that $m = j^2$ and the definition of \hr{j},
        we get the desired.

        Also remember that $r' = -r'_0$ and
        $\sig{p<m}{s_p} = 0$ (see ~\ref{clm_sum}). We get:

        \[
        \begin{split}
            & \norm{r_j}{0} = \norm{-r'_0 + \hr{1} - \hr{1} +
            \ldots + \hr{j-1} - \hr{j-1}
            + \hr{j} + \sum_{p<m}s_p}{0} =
            \\
            & \norm{-\hr{0} + \hr{1} - \hr{1} +
            \ldots + \hr{j-1} - \hr{j-1}
            + \hr{j} + \sum_{p<m}s_p}{0} \le
            \\
            &   \norm{-\hr{0} + \hr{1} + \sum_{p<j}s_p}{0} +
                \norm{-\hr{1} + \hr{2} + \sum_{j \le p<2j}s_p}{0} +
            \\
            &   + \ldots +
                \norm{-\hr{j-1} + \hr{j} + \sum_{j(j-1) \le p<m}s_p}{0} =
            \\
            &   \norm{h_1(-\hr{0} + \hr{1} +
                    \sum_{p<j}s_p)}{\Bm{0}{1}} +
                \norm{h_1(-\hr{1} + \hr{2} +
                    \sum_{j \le p<2j}s_p)}{\Bm{1}{2}} +
            \\
            &   + \ldots +
                \norm{h_1(-\hr{j-1} + \hr{j} +
                    \sum_{j(j-1) \le p<m}s_p)}{\Bm{j-1}{j}} =
            \\
            &   \norm{h_1(-\rt{0} + \rt{j} +
                    \sum_{p<j}s_p)}{\Bm{0}{j}} +
                \norm{h_1(-\rt{j} + \rt{2j} +
                    \sum_{j \le p<2j}s_p)}{\Bm{j}{2j}} +
            \\
            &   + \ldots +
                \norm{h_1(-\rt{(j-1)j} + \rt{m} +
                    \sum_{j(j-1) \le p<m}s_p)}{\Bm{(j-1)j}{m}}
            \\
            \end{split}
        \]

        The last equality is true just by definition of \hr{\l}
        and indiscernibility of \a{i} in \CB.

        Now (remembering that $j>1$)
        we can apply the main claim (\ref{main_clm}) and get
        for each $0 \le \l \le j$ the following inequality:
        \[
            \norm{h_1(-\rt{\l j} + \rt{(\l+1)j} +
            \sum_{\l j \le p < (\l+1)j}s_p)}{\Bm{\l j}{(\l+1)j}}
            \le
            (1+\frac{2}{j})\norm{-\rt{\l j} + \rt{(\l+1)j} +
            \sum_{\l j \le p < (\l+1)j}s_p}{1}
        \]

        Therefore,

        \[
        \begin{split}
            & \norm{r_j}{0} \le
                (1+\frac{2}{j})
                (\norm{-\rt{0} + \rt{j} +
                    \sum_{p<j}s_p}{1} +
                \norm{-\rt{j} + \rt{2j} +
                    \sum_{j \le p<2j}s_p}{1} +
            \\
            &   + \ldots +
                \norm{-\rt{(j-1)j} + \rt{m} +
                    \sum_{j(j-1) \le p<m}s_p}{1})
        \end{split}
        \]

        Rewriting the last inequality in a different way, we get
        \[
        \begin{split}
            & (1+\frac{2}{j})^{-1}\norm{r_j}{0} \le
            \\
            &   \norm{-\rt{0} + \rt{j} +
                    \sum_{p<j}s_p}{1} +
                \norm{-\rt{j} + \rt{2j} +
                    \sum_{j \le p<2j}s_p}{1} +
            \\
            &   + \ldots +
                \norm{-\rt{(j-1)j} + \rt{m} +
                    \sum_{j(j-1) \le p<m}s_p}{1} \le
            \\
            &   \norm{-\rt{0} + \rt{1} + s_0}{1} +
                \norm{-\rt{1} + \rt{2} + s_1}{1} + \ldots +
                \norm{-\rt{m-1} + \rt{m} + s_{m-1}}{1} =
            \\
            &   \norm{r_m}{1}
        \end{split}
        \]

        Which finishes the proof of ~\ref{clm_conv} (2) for the
        case $i=1$. A similar argument is used for $i=-1$, and we
        are done.

    \end{enumerate}
    \end{prf}

    By ~\ref{clm_conv} (1), each one of the three sequences
    \seq{\norm{r_m}{i} : m < \om}  converges. By ~\ref{clm_conv} (2), all of them
    converge to the same limit. Let us denote this limit by $\rho(r',r'') \in \BR$.

    Let $V$ be an ultraproduct of all the $V_m$ modulo some nonprincipal
    ultrafilter $\FD$ on $\om$ (where $V_m$ is a normed space with the norm
    \norm{\cdot}{1}):
    \[
        V = \prod_{m<\om} V_m / \FD
    \]
    \begin{clm_rmk}
    \begin{enumerate}
    \item
        Certainly, this is where the compactness becomes important. We will use
        several times the analog of \L o\`s's theorem for positive bounded formulae,
        claiming $V \models \ph(\infseq{\xi})$ if and only if $V_i \models
        \ph(\xi_i)$ for ``almost all'' $i$.
    \item
        Instead of looking at $V_m$ and $V$, we should have looked at
        their completions, which are Banach spaces, and not just normed
        spaces, but it doesn't matter.
    \item
        We will think of $V$ as embedded into our ``monster'' $\CB$.
    \item
        Note that there is a natural embedding $i_m$ of $V_m$ into $V$:
        \[
            i_m(r) = (0, \ldots, 0, r, \ldots, r, \ldots)
        \]
        i.e. $i_m(r) = \func{g}{\om}{\bigcup V_m} \;$ s.t.
        $\; g(k) = 0 \tfor k<m \tand g(k) = r \tfor k \ge m$.
        Moreover, for $k<m$ we get $i_m \upharpoonright k = i_k$.

        So we will not distinguish between elements of $V_m$ for some $m$
        (in fact, for all $k \ge m$) and the appropriate elements of $V$.
    \end{enumerate}
    \end{clm_rmk}
    The following discussion will be done inside $V$ (and therefore inside $\CB$).
    Let $\b{\om} \in V$ be the ``limit'' of the sequence \seq{\b{m} : m \in \om}, i.e.
    $\b{\om} = \seq{\b{m} : m \in \om} / \FD$. Let $g_\om$ be the ``limit'' of
    \seq{g_m: m \in \om} taking $\b{0}$ onto $\b{\om}$.

    \begin{clm}
    \label{clm_sym}
        Let
        $r', r'' \in \bgen{0}{}$, define
        $r_1 = r' + g_\om(r'')$, $r_{-1} = r'' + g_\om(r')$.
        Let $\rho = \rho(r', r'')$.
        Then $\norm{r_1}{V} = \norm{r_{-1}}{V} = \rho$.
    \end{clm}
    \begin{prf}
        Denote $r_m = r' + g_m(r'')$, $r_{-m} = r'' + g_m{r'}$. Then
        \[
            \norm{r_m}{V_m} = \norm{r_m}{1},\;\;
            \norm{r_{-m}}{V_m} = \norm{r_m}{-1}
        \]
        Remember that by ~\ref{clm_conv}, both \norm{r_m}{1} and \norm{r_m}{-1}
        are ascending sequences converging to \rh. So on one hand,
        $\norm{r_m}{1} \le \rho \tand \norm{r_m}{-1} \le \rho \;$
        for all $m$, and therefore
        \begin{equ}
        \label{equ_less_rh}
        \[
            \norm{r_1}{V} \le \rho, \;\;
            \norm{r_{-1}}{V} \le \rho
        \]
        \end{equ}
        On the other hand, for every real $\eps > 0$, for almost all $m$,
        \[
            \norm{r_m}{1} \ge \rh-\eps, \;\;
            \norm{r_m}{-1} \ge \rh-\eps
        \]
        Therefore, for all real $\eps>0$,
        \begin{equ}
        \label{equ_greater_rh}
        \[
            \norm{r_1}{V} \ge \rho-\eps, \;\;
            \norm{r_{-1}}{V} \ge \rho-\eps
        \]
        \end{equ}
        Combining \ref{equ_less_rh} with \ref{equ_greater_rh}, we get the desired
        equalities.
    \end{prf}

    The following claim can be viewed as the heart of the proof
    we've been working hard for:
    \begin{imp_clm}
    \label{clm_imp}
        $\tp(\b{0},\b{\om})$ is symmetric, i.e. $\tp(\b{0},\b{\om}) = \tp(\b{\om},\b{0})$
    \end{imp_clm}
    \begin{prf}
        Define the obvious mapping $\Phi$ from \dbgen{0}{\om}{} onto itself, extending
        $g_\om \cup g_\om^{-1}$ (``exchanging'' $\b{0} \tand \b{\om}$ and respecting the
        linear structure). It is obviously an isomorphism of vector spaces, so we
        just have to show it is also an isometry. Take $r \in \dbgen{0}{\om}{}$,
        then for some $r', r'' \in \bgen{0}{}$, $r = r' + g_\om(r'')$. Therefore
        $\Phi(r) = r'' + g_\om(r')$. Now, by ~\ref{clm_sym},
        $\norm{r}{V} = \norm{\Phi(r)}{V} = \rh(r', r'')$, q.e.d.
    \end{prf}
    Now we have obviously reached a contradiction. Why? First of all,note that as
    $\tp(\b{i},\b{i+1}) = \tp(\a{i},\a{i+1})$ for all $i \in \om$, we get
    $q(\b{i},\b{i+1})$ for all $i$ (remember: $q(\bar{x},\bar{y})$ defines a partial
    order on \CB, \infseq{\bar a} is ordered by $q$). As $q(\bar{x},\bar{y})$ is a
    partial order, it is in particular transitive, so $q(\b{0},\b{m})$ holds for all
    $m>0$, and therefore $\CB \models q(\b{0},\b{\om})$ (compactness + the fact that
    $q$ is a positive bounded type). But by ~\ref{clm_imp},$\CB \models q(\b{\om},\b{0})$, a
    contradiction to $q$ being a partial order!

\end{prf}

Note that the only property of positive bounded formulae we used in the proof is
that they satisfy the compactness theorem, therefore in fact we proved:

\begin{thm}
    Let \CL be a logic satisfying the compactness theorem for the AEC
    of Banach spaces. Then the class of Banach spaces does not have
    the \sop{} in \CL.
\end{thm}. \qed

\section{Groups}
Let $\CG$ be the ``monster'' group (the universal domain). Our first theorem in this
section is a non-structure result that once again doesn't seem to be surprising as
after seeing the undecidability of the word problem, we feel that any ``bad''
syntactic property can be somehow found in the class of groups.

\begin{prp}
    $\CG$ has $SOP_3$
\end{prp}
\begin{prf} Consider the formula \dphi{x}{y} defined by
``$(\comut{x}{y}) \land (x \neq y)$''.
\begin{itemize}
\item[*] First, we have to show that there is a sequence
\infseq{a} such that $i<j \then \dphi{a_i}{a_j}$. But this is trivial by using HNN
extentions and compactness.
\item[*] Secondly we have to make sure there is no ``triangle'', but
this is actually a well-known example in geometric group theory (see
[Grp]) of a triangle $X = \langle a,b: \comut{a}{b} \rangle,
Y = \langle b,c: \comut{b}{c} \rangle, Z = \langle a,c: \comut{c}{a} \rangle$
that generates a \emph{trivial} group when put together. Therefore, \\
\mbox{$\CG \models (\forall x, y, z)(\comut{x}{y} \land \comut{y}{z} \land
\comut{z}{x} \longrightarrow x = y = z = e)$}
(where $e$ is the group identity).
Therefore \mbox{$\CG \models \neg (\exists x, y, z)( \dphi{x}{y} \land
\dphi{y}{z} \land \dphi{z}{x}) $}, as required.
\end{itemize}
\end{prf}

The proof uses the fact that there can not be a triangle of a certain kind. A
natural question now is - what about quadriangles? In particular, is the group $H =
\langle a,b,c,d : \comut{a}{b}, \comut{b}{c}, \comut{c}{d}, \comut{d}{a} \rangle$
also trivial? Once again, it's a well-known fact that it is actually infinite, and
the proof is even more interesting than the fact itself, as it seems very general -
it doesn't speak at all about the relations between the generators. In fact, the
proof suggests a generalization that roughly speaking says that it is impossible to
``collapse'' a group with four generators by forcing relations between only adjacent
pairs. Model theoretically, it leads to the following surprising \emph{structure}
result, showing that unlike what people might have thought, there is a hope for some
model-theoretic structure theory for the class of all groups.

\begin{thm}
\CG does \emph{not} have $SOP_4$
\end{thm}
\begin{prf}
Suppose towards a contradiction that \dbarphi{x}{y} exemplifies \sop{4} in \CG. In
particular, there exists an indiscernible sequence \infseq{\bar a} such that $i<j
\then \dbarphi{a_i}{a_j}$. Define for all $i \in \omega$, $H_i = \gen{\bar
a_i}{\CG}$. We denote $\len(\bar x) = \len(\bar y)$ in \dbarphi{x}{y} by $n$ and
assume wlog (by indiscernibility) that there exists $n^* < n$ such that
$\bigwedge_{\ell<n^*} (a_{i,\ell} = a^*_\ell)$ for all $i < \om$ and \seq{
a_{i,\ell} : n^* \le \ell < n, i < \om} is a sequence of distinct elements. Define
$H^- = \seq{\subsup{a}{\l}{*}: \l < n^*}$, i.e. $H^- = H_i \cap H_j$ for all
$i<j<\om$.

By the indiscernibility , there exists for $i \neq j \in \omega$, an isomorphism
\func{f_{i,j}}{H_i}{H_j} mapping $\bar a_i$ onto $\bar a_j$. Define for all $i<j \in
\omega$, $H_{i,j} = \gen{\bar a_i, \bar a_j}{\CG}$. For $j<i \in \omega$ we define
\dH{i}{j} by ``relabeling'', changing the roles of $\bar a_i$ and $\bar a_j$, i.e.
as a set \dH{i}{j} equals \dH{j}{i}, and the group action is defined on it such that
there exists \func{f_{j,i}^{i,j}}{\dH{j}{i}}{\dH{i}{j}} an isomorphism extending
$f_{i,j} \cup f_{j,i}$. So for $j<i$, \dH{i}{j} does not have to be a subgroup of
\CG (but we can embed it into \CG,  as  \CG  is universal).

Given two groups $G_1$ and $G_2$ and a subgroup of both, $G_0$, we shall denote the
free amalgamation of the two over $G_0$ by \mbox{$\famalg{G_1}{G_2}{G_0}$}. Now let
us concentrate on $H_0, H_1, H_2, H_3$. Define \mbox{$K_0 =
\famalg{H_0}{H_2}{H^-}$}, \mbox{$K_1 = \famalg{H_{0,1}}{H_{1,2}}{H_2,H^-}$},
\mbox{$K_2 = \famalg{H_{2,3}}{H_{3,0}}{H_3,H^-}$}. Once again, those groups do not
have to be subgroups of \CG. It is obvious, though, that $K_0$ is a subgroup of both
$K_1$ and $K_2$ (by definition of free product and amalgamation of groups). So we
define $K = \famalg{K_1}{K_2}{K_0}$.

\CG is universal, so we can embed $K$ into \CG. Denote the image of
$\a{i}$ under this embedding by $\bar b_i \in \CG$. Now we note
\begin{clm}
\label{clm:grptp} $\tp(\b{0} \b{1}, \CG) = \tp(\b{1} \b{2}, \CG) =  \tp(\b{2} \b{3},
\CG) = \tp(\b{3} \b{0}, \CG) = \tp(\a{0} \a{1}, \CG)$.
\end{clm}
\begin{prf}
$\tp(\b{0} \b{1}, \CG) = \tp(\a{0} \a{1}, K) = \tp(\a{0} \a{1}, K_1) = \tp(\a{0}
\a{1}, H_{0,1}) = \tp(\a{0} \a{1}, \CG)$. The first equality is true because types
are preserved under group isomorphisms (``embeddings''), and the rest - just the
definitions of the groups. Using the same arguments for \dtpin{\b{1}}{\b{2}}{\CG},
we get $\dtpin{\b{1}}{\b{2}}{\CG} = \dtpin{\a{1}}{\a{2}}{\CG}$, but the latter
equals \dtpin{\a{0}}{\a{1}}{\CG} by indiscernibility. The same argument (replacing
$K_1$ by $K_2$) shows $\dtpin{\b{2}}{\b{3}}{\CG} = \dtpin{\a{0}}{\a{1}}{\CG}$. Now
$\tp(\b{3} \b{0}, \CG) = \tp(\a{3} \a{0}, K) = \tp(\a{3} \a{0}, K_2) = \tp(\a{3}
\a{0}, H_{3,0})$, but the latter equals
$\dtpin{f_{0,3}^{3,0}(\a{3})}{f_{0,3}^{3,0}(\a{1})}{\CG} = \tp(\a{0} \a{3}, \CG)$ by
the definition of $H_{3,0}$ and $f_{0,3}^{3,0}$, and by indiscernibility we're done.
\end{prf}
Now we obviously get a contradiction, as by (\ref{clm:grptp})
$\CG \models \dbarphi{b_0}{b_1} \land   \dbarphi{b_1}{b_2} \land
\dbarphi{b_2}{b_3} \land \dbarphi{b_3}{b_0}$, which contradicts the fact
that \dbarphi{x}{y} exemplifies \sop{4} in \CG.
\end{prf}

\end{document}